\documentclass[conference, a4paper,UTF8]{IEEEtran}
\IEEEoverridecommandlockouts
\usepackage{fancyhdr} 
\usepackage{color}
\usepackage{amssymb}
\usepackage{amsmath}
\usepackage{amsmath,amssymb,amsfonts}
\usepackage{algorithmic}
\usepackage{graphicx}
\usepackage{textcomp}
\usepackage{amsmath,amssymb,amsfonts}
\usepackage{algorithmic}
\usepackage{graphicx}
\usepackage{float}
\usepackage{subfigure}
\usepackage{setspace}
\usepackage{hyperref}
\usepackage{textcomp}
\usepackage{amsfonts}
\usepackage{paralist}
\usepackage{tikz}
\usetikzlibrary{shapes.arrows}
\usepackage{amsmath}
\usepackage{algorithm}
\usepackage{algorithmic}
\usepackage{amsthm}
\usepackage{warpcol}
\usepackage{booktabs}
\usepackage{multirow}
\usepackage{graphicx}
\usepackage{booktabs}
\usepackage{rotating}
\usepackage{placeins}
\usepackage{mathrsfs}
\usepackage{ifsym}
\usepackage{tikz} 
\usepackage{subfigure}
\usepackage{listings}

\usepackage{geometry}
\usetikzlibrary{decorations.pathreplacing}
\usetikzlibrary{decorations.pathreplacing,calligraphy}
\usepackage{makecell}
            
\usepackage{lastpage}      
\usepackage{setspace}          
\usepackage{algorithm} 
\usepackage{algorithmic} 
\usepackage{multirow} 

\usepackage{amsmath,bm} 
\usepackage{xcolor}                         
\usepackage{layout}  
\usepackage{amsmath}   
\usepackage{hyperref}
\usepackage{color}
\usepackage{amsfonts}
\usepackage{graphicx}       
\usepackage{hyperref}                                  
\begin{document}
\begin{spacing}{1}

\title{Solving the Three-Player-Game}
\author{
\IEEEauthorblockN{Fangqi Li}
\IEEEauthorblockA{\textit{School of Cyber Science and Engineering, SEIEE, SJTU}\\
\IEEEauthorblockN{$\{$solour\_lfq$\}$@sjtu.edu.cn}}
}
\date{Today}
\maketitle
\begin{abstract}
In this paper we solve the three-player-game question. A three-player-game consists of a series of rounds. There are altogether three players. Two players participate in each round, at the end of the round the loser quits and the third player enters the ring and another round starts. The game terminates if all six win-lose relationships appear. During each round, two players win with equal probability. One is asked to calculate the expectation of the number of rounds. It turns out to be an exemplary question that involves probabiltiy theory and dynamic programming. It can serve as an instance or exercise in the chapter of conditional expectation of any elementary or advanced textbook on probability.


\end{abstract}
\section{Introduction}
\label{section:1}
The three-player-game (3PG) is an interesting mathematical quiz. A 3PG involves three players and consists of a series of rounds, each zero-win round involves two out of the three players, at the end of the round the loser quits and the third player enters the ring and another round starts. The game terminates if all six win-lose relationships appear. During each round, two players win with equal probability. 

For example, let Alice, Bob and Carole be three players. Alice and Bob play the first round, then Alice loses. Therefore the second round involves Bob and Carole. If Carole loses then Alice has to return to the ring. The game proceeds until every player has beaten the other two players for at least once. 

One should take the win-lose situation of each round as the underlying infinite probability space and the number of rounds as a random variable. The target is to compute the expectation of this random variable. 

We solve this question from scratch and give some additional analysis on generalization.

\section{Solving 3PG by Reduction}
\label{section:2}
\subsection{Formulation and Reduction Rules}
It is straightforward to observe that the sufficient statistics of any round in the game is the occurance of all possible win-lose relationships sofar and the current players in the ring. That is to say, at each stage of the game, the current situation is comprehensively described by six bits (indicate whether or not the six win-lose relationship appear) and one ternary bit (indicate the current players in the ring). We denote all possible states in the following form:
$$\mathcal{S}=\left\{\left(b_{1}b_{2}b_{3}b_{4}b_{5}b_{6}|t_{1}t_{2}\right)\right\},$$
where each $b_{i}\in\left\{0,1\right\}$ denotes whether a win-lose relationship appears or not. Whereas $t_{1}$ and $t_{2}$ in $\left\{1,2,3 \right\}$ are the indices for the players in battle. From $i=1,\cdots,6$, $b_{i}$ denotes the win-lose relationship of:
$$\text{Player 1 beats Player 2},$$
$$\text{Player 2 beats Player 1},$$
$$\text{Player 1 beats Player 3},$$
$$\text{Player 3 beats Player 1},$$
$$\text{Player 2 beats Player 3},$$
$$\text{Player 3 beats Player 2}.$$
For each $S\in\mathcal{S}$, let $f(S)$ be the expectation of the number of rounds of a modified 3PG begins with $S$. Now we are asked to compute, w.l.o.g.,:
$$f(000000|\left\{1,2\right\}).$$

It is natural to use conditional expectation to introduce the reduction rule, we begin with the definition:
$$\mathbb{E}[X]=\mathcal{P}(A)\mathbb{E}[X|A]+\mathcal{P}(\overline{A})\mathbb{E}[X|\overline{A}],$$ 
where $X$ is any r.v. and $A$ is an event. Let $X$ be \emph{the number of rounds of a 3PG begins with $S$} and $A$ be the event that $t_{1}$ wins over $t_{2}$ w.r.t. $S$, we have:
$$\mathbb{E}[X]=f(S),$$
$$\mathcal{P}(A)=\mathcal{P}(\overline{A})=\frac{1}{2},$$
$$\mathbb{E}[X|A]=1+f(S_{1}),$$
$$\mathbb{E}[X|\overline{A}]=1+f(S_{2}),$$
where $S$ reduces to state $S_{1}(S_{2})$ if $A$ does (not) happen. This expansion yields:
\begin{equation}
\label{equation:1}
f(S)=1+\frac{1}{2}\cdot f(S_{1})+\frac{1}{2}\cdot f(S_{2}),
\end{equation}
where $S_{1}$ and $S_{2}$ are two results from $S$ given $t_{1}$ or $t_{2}$ wins. \eqref{equation:1} would become the basis of reduction which finally solves the 3PG. We are going to see a dozen of examples in the coming section. 

The reduction terminates at the basic states
$$f(111111|\left\{t_{1},t_{2}\right\})=0,$$
where $t_{1},t_{2}$ can be any players.

We define the \emph{order} of a state by the number of 0s in its binary parts. For example, $(000000|\left\{1,2\right\})$ is a state of order six. 

The number of states $|S|$ is 192, however, some states is not going to appear with the root state $(000000|\left\{1,2\right\})$. Moreover, one is encouraged to evoke symmetry to further reduce computation. We define two states $S_{A}$ and $S_{B}$ to be symmetric if there exist a permutation $\pi$ on $\left\{1,2,3\right\}$ such that 
$$\pi(S_{A})=S_{B},$$
where applying $\pi$ on a state $S\in\mathcal{S}$ changes both the order of binary indicators and the name of players in the ring. For $b$ indicates the relationship Player 1 beats Player 2 in $S_{A}$, its values also indicates whether the relationship Player $\pi (1)$ beats Player $\pi (2)$ in $S_{B}$. That is to say, $\pi$ on $\left\{1,2,3\right\}$ introduces a permutation on $\left\{1,2,3,4,5,6 \right\}$, for example, $\pi =(1,2)$ introduces (using group algebra notation):
$$(1,2)(3,5)(4,6)$$
While $\pi=(2,3)(1,2)$ introduces:
$$(1,4,5)(2,3,6)$$
on the first six binary bits on $\mathcal{S}$. Finally, 
$$t_{B,i}=\pi (t_{A,i}),i=1,2.$$

For example, let
$$S_{A}=(101000|\left\{1,3 \right\})$$
$$S_{B}=(101000|\left\{1,2 \right\})$$
then $S_{A}$ and $S_{B}$ are symmetric by adopting
$$\pi=(2,3).$$
Naturally, symmetry is an equivalent relationship, and for symmetric states $S_{A}$ and $S_{B}$ we have:
$$f(S_{A})=f(S_{B}),$$
since their difference is only a matter of naming. This observation helps to reduce the number of states significantly. However, there is hardly any method to examine whether two states are symmetric other than checking all possible permutations. 

\subsection{Preparations and Preprocessing}
Having obtained \eqref{equation:1}, one might eagerly argue that a simple recursion program would trivially solve the task:

Define:
$$\text{Compute}(S)=1+\frac{\text{Compute}(S_{1})}{2}+\frac{\text{Compute}(S_{1})}{2},$$
$$\text{Compute}(111111|\left\{t_{1},t_{2} \right\})=0.$$
Return $\text{Compute}(000000|\left\{1,2\right\})$.

With $\text{Compute}(\cdot)$ as an algorithmatic realization of $f(\cdot)$. Computing $\text{Compute}(000000|\left\{1,2\right\})$ is then done automatically by spanning a recursion tree (let $S_{1},S_{2}$ be two children of the node represents $S$), during which dynamic programming might help to reduce computation time \cite{cormen2009introduction}. However, this method is not determined to success since \eqref{equation:1} does not ensured that a state $S$ itself does not appear in the computing tree spanned by with $S$ as the root, which fact deadlocks this paradigm. At certain states, it is necessary to use the linear relationship between their expectations to solve $f(\cdot)$ and a naive recursion is far from enough. 

This fact, together with the last observation from the previous section, indicates that instead of \emph{passively} spanning a computing tree, we should better \emph{aggressively} compute the leave states (those states with few 0s in their binary part) at first. 

Before actually conducting reduction from $(000000|\left\{1,2\right\})$, we conduct preprocessing by computing $f(S)$ for some elementary states (states with small orders) beforehand. These computations are collected into a series of gradual propositions. 

\textbf{Proposition A:} 
$$f(111110|\left\{2,3 \right\})=4,$$
$$f(111110|\left\{1,2 \right\})=6,$$
$$f(111110|\left\{1,3 \right\})=6,$$

\textbf{Proof:} Let $x,y,z$ denote these three values respectively, we have according to \eqref{equation:1}:
$$x=1+\frac{y}{2},$$
$$y=1+\frac{z}{2}+\frac{x}{2},$$
$$z=1+\frac{y}{2}+\frac{x}{2}.$$
That is to say 
$$
\begin{pmatrix}
x \\
y\\
z
\end{pmatrix}=
\begin{pmatrix}
0 & \frac{1}{2} & 0\\
\frac{1}{2} & 0 & \frac{1}{2}\\
\frac{1}{2} & \frac{1}{2} & 0\\
\end{pmatrix}
\begin{pmatrix}
x \\
y\\
z
\end{pmatrix}+
\begin{pmatrix}
1 \\
1\\
1
\end{pmatrix}
$$
This gives $x=4,y=6,z=6$ as the only solution. \qed

Proposition A finishes the computation of all states $S$ with five 1s in their binary part, i.e., all states of order one. Technically, let the only component as 0 be $t_{A}$ beats $t_{B}$ in a state $S=(\cdots|\left\{t_{1},t_{2}\right\})$, then if $\left\{t_{1},t_{2} \right\}=\left\{t_{A},t_{B} \right\}$ then the expectation of the corresponding state $f(S)$ is 4, otherwise it is 6. We now use this as the block of building estimation for states with order two. 

Let us begin with states with two vacant relationships $t_{A}$ beats $t_{B}$ and $t_{B}$ beats $t_{A}$.

\textbf{Proposition B:}
$$f(111100|\left\{2,3 \right\})=7.$$

\textbf{Proof:} Applying \eqref{equation:1} to this state:
$$
\begin{aligned}
f(111100|\left\{2,3\right\})=1&+\frac{1}{2}\cdot f(111110|\left\{1,2\right\})\\
&+\frac{1}{2}\cdot f(111101|\left\{1,3\right\}).
\end{aligned}
$$
Applying Proposition A finishes the proof. \qed

\textbf{Proposition C:}
$$f(111100|\left\{1,3\right\})=9.$$

\textbf{Proof:}
Applying \eqref{equation:1} to this state:
$$
\begin{aligned}
f(111100|\left\{1,3\right\})=1&+\frac{1}{2}\cdot f(111100|\left\{1,2\right\})\\
&+\frac{1}{2}\cdot f(111100|\left\{2,3\right\}).
\end{aligned}
$$
Now considering $\pi =(2,3)$ then we have:
$$f(111100|\left\{1,3\right\})=f(111100|\left\{1,2\right\}),$$
combining this with Proposition B finishes the proof.\qed

Now if the only two vacant relationships are $t_{A}$ beats $t_{B}$ and $t_{B}$ beats $t_{A}$, we are ready to read the expectation of the state $S(\cdots|\left\{t_{1},t_{2}\right\})$. If $\left\{t_{1},t_{2}\right\}=\left\{t_{A},t_{B}\right\}$ then $f(S)$ is 7, else it is 9.

We then proceed to states with vacant relationships 
\begin{itemize}
\item $t_{A}$ beats $t_{B}$ and $t_{A}$ beats $t_{C}$, 
\item $t_{A}$ beats $t_{B}$ and $t_{C}$ beats $t_{B}$, 
\item $t_{A}$ beats $t_{B}$ and $t_{B}$ beats $t_{C}$. 
\end{itemize}

\textbf{Proposition D:}
$$f(010111|\left\{2,3 \right\})=8,$$
$$f(010111|\left\{1,2 \right\})=7,$$

\textbf{Proof:} Let $x,y$ denote $f(010111|\left\{2,3 \right\})$ and $f(010111|\left\{1,2 \right\})$, according to \eqref{equation:1} (one easily notes that $(010111|\left\{1,2 \right\})$ is symmetrical to $(010111|\left\{1,3 \right\})$):
$$x=1+\frac{y}{2}+\frac{y}{2},$$
$$y=1+\frac{f(110111|\left\{1,3 \right\})}{2}+\frac{x}{2}.$$ Applying Proposition A finishes the proof. \qed

\textbf{Proposition E:}
$$f(101011|\left\{2,3 \right\})=9,$$
$$f(101011|\left\{1,2 \right\})=8,$$

\textbf{Proof:} Let $x,y$ denote $f(101011|\left\{2,3 \right\})$ and $f(101011|\left\{1,2 \right\})$, using symmetry and Proposition A as in the proof of Proposition D:
$$x=1+y,$$
$$y=1+\frac{y}{2}+\frac{6}{2}.$$
This finishes the proof.\qed

\textbf{Proposition F:}
$$f(011101|\left\{1,2\right\})=\frac{38}{5},$$
$$f(011101|\left\{2,3\right\})=\frac{36}{5},$$
$$f(011101|\left\{1,3\right\})=\frac{42}{5}.$$
\textbf{Proof:} Let $x,y,z$ denote these three values, using \eqref{equation:1} and Proposition A:
$$x=1+\frac{6}{2}+\frac{y}{2},$$
$$y=1+\frac{4}{2}+\frac{z}{2},$$
$$z=1+\frac{x}{2}+\frac{y}{2}.$$
Solving this system yields
$$x=\frac{38}{5},y=\frac{36}{5},z=\frac{42}{5}.$$
\qed

Proposition B-F finish computing states of order two.  

The number of states with three 1s/of order three is larger. There are at least $\frac{3*\binom{6}{3}}{3\!}=10$ unsymmetric states. Although unnecessary for the following sections, one is encouraged to compute all 13 independent states with three appeard win-loss relationships. 

We are now ready to begin from $f(000000|\left\{1,2\right\})$ and hope that compution meets with the Propositions A-F before at an early stage of computation. 

\subsection{The Main Reduction}
Attempting to solve this problem by reduction, applying \eqref{equation:1} onto $f(000000|\left\{1,2 \right\}).$:
$$
\begin{aligned}
f(000000|\left\{1,2\right\})=1&+\frac{1}{2}\cdot f(100000|\left\{1,3\right\})\\
&+\frac{1}{2}\cdot f(010000|\left\{2,3 \right\}).
\end{aligned}
$$

From now on, let $S_{1}$ be the state where the player with the smaller index winning the current round. Adopting
$$\pi =(1,2),$$
then $(100000|\left\{1,3\right\})$ and $(010000|\left\{2,3 \right\})$ are symmetric, so:
\begin{equation}
\label{equation:2}
f(000000|\left\{1,2\right\})=1+f(100000|\left\{1,3\right\}).
\end{equation}
Therefore we are left with the problem of computing $f(100000|\left\{1,3\right\})$. Now
\begin{equation}
\label{equation:3}
\begin{aligned}
f(100000|\left\{1,3\right\})=1&+\frac{1}{2}\cdot f(101000|\left\{1,2\right\})\\
&+\frac{1}{2}\cdot f(100100|\left\{2,3 \right\}).
\end{aligned}
\end{equation}
We first address $f(101000|\left\{1,2\right\})$ and then return to {\color{red}$f(100100|\left\{2,3\right\})$}. Since
$$
\begin{aligned}
f(101000|\left\{1,2\right\})=1&+\frac{1}{2}\cdot f(101000|\left\{1,3\right\})\\
&+\frac{1}{2}\cdot f(111000|\left\{2,3 \right\}).
\end{aligned}
$$
However, let 
$$\pi =(2,3),$$ 
we conclude that
$$f(101000|\left\{1,2 \right\})=f(101000|\left\{1,3\right\}),$$
thus 
\begin{equation}
\label{equation:4}
f(101000|\left\{1,2 \right\})=2+f(111000|\left\{2,3\right\}).
\end{equation}

Keep reducing:
$$
\begin{aligned}
f(111000|\left\{2,3\right\})=1&+\frac{1}{2}\cdot f(111010|\left\{1,2\right\})\\
&+\frac{1}{2}\cdot f(111001|\left\{1,3 \right\}).
\end{aligned}
$$
We are now meeting two states with four 1s and two 0s, $f(111010|\left\{1,2\right\})$ is symmetric to $f(010111|\left\{2,3\right\})$, hence its value is 8 according to Proposition D, while $f(111001|\left\{1,3\right\})$ addresses a state symmetric to $f(011101|\left\{2,3 \right\})$, whose value is $\frac{36}{5}$ according to Proposition F. Pluggin them back into \eqref{equation:4} gives:
$$f(101000|\left\{1,2\right\})=\frac{53}{5}.$$

To return to \eqref{equation:3}, we still need to compute $f(100100|\left\{2,3\right\})$.
\begin{equation}
\label{equation:5}
\begin{aligned}
f(100100|\left\{2,3\right\})=1&+\frac{1}{2}\cdot f(100110|\left\{1,2 \right\})\\
&+\frac{1}{2}\cdot f(100101|\left\{ 1,3\right\}).
\end{aligned}
\end{equation}
For 
$$
\begin{aligned}
f(100110|\left\{1,2\right\})=1&+\frac{1}{2}\cdot f(100110|\left\{1,3 \right\})\\
&+\frac{1}{2}\cdot f(110110|\left\{2,3 \right\}).
\end{aligned}
$$

Where as one can easily observe the symmetry between $(100110|\left\{1,2\right\})$ and $(100110|\left\{1,3\right\})$, we have:
$$f(100110|\left\{1,2\right\})=2+f(110110|\left\{2,3\right\})=\frac{46}{5}.$$
The last term remained is $f(100101|\left\{1,3\right\})$. We begin with
$$
\begin{aligned}
f(100101|\left\{1,3\right\})=1&+\frac{1}{2}\cdot f(101101|\left\{1,2\right\})\\
&+\frac{1}{2}\cdot f(100101|\left\{2,3 \right\})\\
=\frac{9}{2}&+\frac{1}{2}\cdot f(100101|\left\{2,3\right\}).
\end{aligned}
$$
Finally, we have:
$$f(100101|\left\{2,3\right\})=\frac{24}{5}+\frac{1}{2}\cdot f(100101|\left\{1,3\right\}).$$
This gives:
$$f(100101|\left\{1,3\right\})=\frac{46}{5}.$$
Pluggin them into \eqref{equation:5} yields:
$$f(100100|\left\{2,3\right\})=\frac{51}{5}.$$
Now \eqref{equation:3} yields:
$$f(100000|\left\{1,3\right\})=11.4.$$
At length, pluggin this into \eqref{equation:2} yields:
$$f(000000|\left\{1,2\right\})=12.4.$$

Hitherto we have finished all the reduction. 

\subsection{Simulation Results}
\begin{figure}[htb]
\centering
\includegraphics[width=0.5\textwidth] {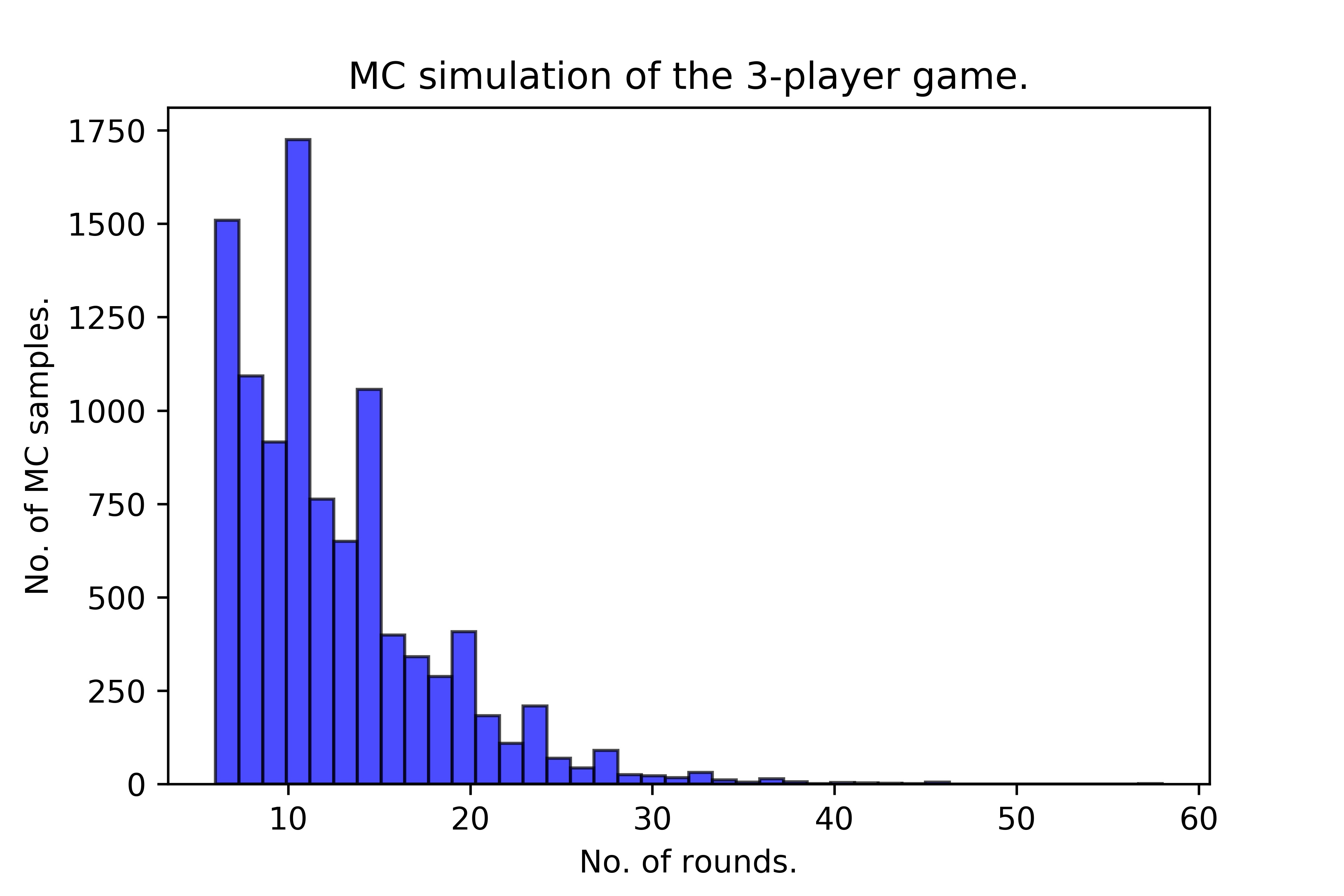}
\caption{Monte Carlo simulation for 3PG with 1,000 samples.} 
\label{figure:1}
\end{figure}

As for an empirical verification of the result, a straightforward Monte Carlo simulation was conducted (with 1,000 samples) and the result is shown as Figure. \ref{figure:1}, the mean of the number of rounds is 12.4287.

\section{The Second Order Analysis}
Given the expectation of the number of rounds in any states:
$$\left\{f(S):S\in\mathcal{S} \right\},$$
it is straightforward to compute the variance of the number of rounds by reduction. The bridge is:
$$\text{var}[X]=\text{var}[\mathbb{E}[X|Y]]+\mathbb{E}[\text{var}[X|Y]].$$
Now let $X$ be the random variable that denotes the number of rounds in the current state and $Y$ be the indicator of the current competition. Let $g(S)$ be the variance of the number of rounds of a 3PG begins from the state $S$, we have:
$$g(S)=\frac{(f(S_{1})-f(S_{2}))^{2}}{2}+\frac{g(S_{1})+g(S_{2})}{2}.$$
Thus given $\left\{f(S):S\in\mathcal{S} \right\}$ it is straightforward to compute $\left\{g(S):S\in\mathcal{S} \right\}$ (repeat what has been done in the sections before, \emph{reversely computing along the martingle!}) and deduce the variance of the number of rounds in a 3PG.

For example, consider the variances of the number of rounds for states $S_{1}=(111110|\left\{2,3 \right\})$, $S_{2}=(111110|\left\{1,2 \right\})$, $S_{3}=(111110|\left\{1,3 \right\})$. Using Proposition A, we have:
$$\begin{pmatrix}
g(S_{1}) \\
g(S_{2}) \\
g(S_{3})
\end{pmatrix}=
\begin{pmatrix}
0 & \frac{1}{2} & 0 \\
\frac{1}{2} & 0 & \frac{1}{2} \\
\frac{1}{2} & \frac{1}{2} & 0 
\end{pmatrix}
\begin{pmatrix}
g(S_{1}) \\
g(S_{2}) \\
g(S_{3})
\end{pmatrix}
+
\begin{pmatrix}
18 \\
2 \\
2
\end{pmatrix}.
$$
Which yields:
$$\begin{pmatrix}
g(S_{1}) \\
g(S_{2}) \\
g(S_{3})
\end{pmatrix}=
\begin{pmatrix}
40 \\
44 \\
44
\end{pmatrix}.
$$

\section{The Probabilistic Framework}
The analysis so far is hardly relied on the probability space. The reason behind is that it is hard to establish the equivalence between an element in the probability space and the value of the random variable \cite{alon2004probabilistic}. Considering:
$$\Omega=\left\{+,- \right\}^{\infty}.$$
Where $+/-$ denotes the player with larger/smaller index winning the current round. To compute $\text{Pr}(X=n)$, where $X$ is the random variable that counts the number of rounds until termination. One has to find the number of $\left\{+,-\right\}^{n}$ sequences where all win-lose relationship appears until the final round. Although it is efficient to transcript a $\left\{+,-\right\}^{n}$ sequence into win-lose relationship sequences, it is hard to write down (be it exists) a tractable necessary and sufficient condition for $X=n$. 

However, we could use the solution of 3PG to answer questions yielded from a more probabilistic perspective. For example:

\emph{Building up a string s with three characters \{`''a'',''b'',''c''\}, s begins with ''a'' and each character is followed by one different character with equal probability, s terminates until all six pairs appears in the string. What is the expected length of s?} This question is isomorphic to 3PG.

\section{Generalization}
Having finished the analysis of 3PG, we now proceed to a genelization study. The problem is, is it possible to find the asymptotical behavior of the solution to $n$-PG? The generalization of 3PG to $n$-PG is not unique, e.g., each round can still involve two players, and one random player enters the next round instead of the loser, or one can adopt ternary logic to mark the result of battles.  

We study the general $n$-PG with two players participating each round, and a player is randomly (uniformly and independently) chosen to replace the loser of the current round in the ring. 

First we try to address the states of order one, w.l.o.g., let the vacant relationship be Player 1 beats Player 2. There are four independent (unsymmetric) states with $\left\{1,2\right\}$, $\left\{1,3\right\}$, $\left\{2,3\right\}$ and $\left\{3,4 \right\}$ as the current pair of players in the ring (assuming $n \geq 4$). Let $x,y,z,w$ denote the corresponding expectations, then we have:
$$
\begin{pmatrix}
x\\
y\\
z\\
w
\end{pmatrix}=
A_{1}
\begin{pmatrix}
x\\
y\\
z\\
w
\end{pmatrix}+
\begin{pmatrix}
1\\
1\\
1\\
1
\end{pmatrix},
$$
with:
$$A_{1}=
\begin{pmatrix}
0 & 0 & \frac{1}{2} & 0\\
\epsilon & \frac{1}{2} & \epsilon & \frac{1}{2}\\
\epsilon & \epsilon & \frac{1}{2} & \frac{1}{2}\\
0 & \epsilon & \epsilon & 1-2\epsilon
\end{pmatrix},$$
where 
$$\epsilon=\frac{1}{2(n-2)}.$$
Since we have:
$$
\begin{pmatrix}
x\\
y\\
z\\
w
\end{pmatrix}=(I-A_{1})^{-1}\cdot
\begin{pmatrix}
1\\
1\\
1\\
1
\end{pmatrix},
$$
as:
$$(I-A_{1})^{-1}=\sum_{i=0}^{\infty}A_{1}^{i},$$
the only task remained is to track the spectral radius of $A_{1}$ \cite{meyer2000matrix}, the trick here is to apply the Gerschgorin theorem to the last row of $A_{1}$, with yields that the largest eigenvalue of $A_{1}$ (assumed to be real) is no less than:
$$\lambda=1-4\epsilon.$$
Therefore the spectral radius of $(I-A)^{-1}$ is no less than:
$$\frac{1}{1-\lambda}\sim O(n).$$
So is the order of $x,y,z,w$. 

Moving to states of order two is a similar case, let the independent states be $x',y',\cdots,w'$, we have:
$$
\begin{pmatrix}
x'\\
y'\\
\cdots\\
w'
\end{pmatrix}=A_{2}
\begin{pmatrix}
x'\\
y'\\
\cdots\\
w'
\end{pmatrix}+
\begin{pmatrix}
c_{x}\\
c_{y}\\
\cdots\\
c_{w}
\end{pmatrix},
$$
Where elements in $(c_{x},c_{y},\cdots,c_{w})^{\text{T}}$ are constants with value 1 or a multiple of $x,y,z,w$ that has been evaluated before. Hence the order of elements in $(c_{x},c_{y},\cdots,c_{w})^{\text{T}}$ is at most $O(n)$. 

To measure the spectral radius of $A_{2}$, we resort to a similar line of reasoning: let $w$ be the expectation of the state where the current players on the ring is different from those players involved in the vacant win-lose relationships. Then the final row of $A_{2}$ has $1-\epsilon_{2}$ as the last component, where
$$\epsilon_{2}=\frac{4}{n-2}$$
in the most probable case. This yields the fact that the order of the spetral radius of $(I-A_{2})^{-1}$ be $O(n)$, hence the order of the expectations of states of order two turns out to be $O(n^{2})$. 

In general, for states of order $\phi$, let $w_{\phi}$ be the expectation of the state where the current players on the ring are free from those $\phi$ vacant pairs, let $A_{\phi}$ be the transition matrix at that order. There are at most $2\phi$ players involve with the vacant pairs, hence the entry on the right-bottom most side of $A_{\phi}$ is at most:
$$1-\frac{2\phi}{n-2}.$$
That is to say, the spectral radius of $(I-A_{\phi})^{-1}$ is of order:
$$O\left(\frac{n}{\phi}\right).$$
Finally, counting all states of order $\phi=0,1,\cdots,n$, we have the order of the solution of an $n$-PG be:
$$O\left(\frac{n^{n}}{n!}\right)=O\left(\frac{\text{e}^{n}}{\sqrt{n}}\right).$$

Analogously, the variance for general $n$-PG can be approximated using the same framework. The $f(\cdot)$ for states of order $\phi$ is of order:
$$\frac{1}{\phi}\left(\frac{n\text{e}}{\phi} \right)^{\phi}.$$
We have that $g(\cdot)$ for states of order $\phi$ (denoted by $g_{\phi}$)is determined by the larger term in $f^{2}_{\phi}$ and $g_{\phi-1}$, so at least:
$$g_{\phi}\geq O\left(\frac{n}{\phi}\cdot f^{2}_{\phi}\right),$$
which is far less than the order of $f^{2}_{\phi}$, therefore we conjecture that the order of the variance in $n$-PG is:
$$O\left(\frac{\text{e}^{2n}}{n} \right).$$

\begin{figure}[htb]
\centering
\includegraphics[width=0.5\textwidth] {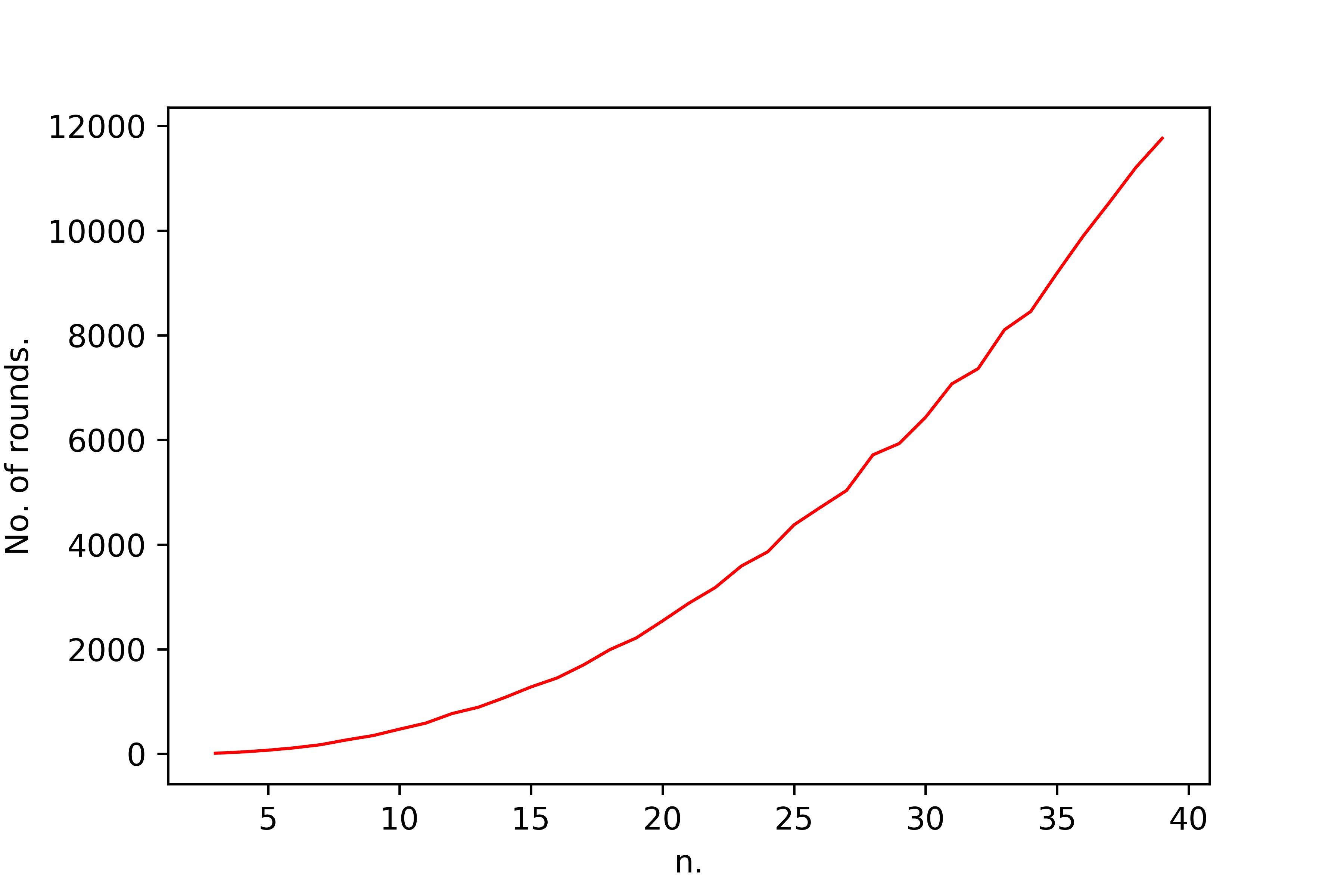}
\caption{Monte Carlo simulation for $n$-PG with 100 samples for each $n$.} 
\label{figure:2}
\end{figure}

In fact, we observe that the estimation on the expectation is possibly a rather slack one, this is due to the following facts:
\begin{itemize}
\item The estimation based on the Gerschgorin theorem on $A$ might significantly increase the spectral radius of $(I-A)^{-1}$.
\item The spectral radius might be involved with negligible terms so the speed of growth declines.
\end{itemize}
The simulation results of general $n$-PGs are illustrated as in Figure. \ref{figure:2}.

From which we might optimistically conjecture that the growth of the expectation of the number of rounds is only of order $n^{2}$, but there seems to be a vacancy in establishing this result. 

\section{Conclusion}
This paper address the 3PG question. We attack this question with dynamic programming, highlight the necessary tricks that signicantly reduce redundant computation and analyze the ideas behind. The general case is also proposed and a rough bound is derived. 

\section{Acknowledgement}
Haoran Ye for provided the 3PG question, Runbo Ni and his colleagues provided an early version of solution. 

\bibliographystyle{ieeetr}
\bibliography{3p.bib}

\begin{thebibliography}{1}

\bibitem{cormen2009introduction}
T.~H. Cormen, C.~E. Leiserson, R.~L. Rivest, and C.~Stein, {\em Introduction to
  algorithms}.
\newblock MIT press, 2009.

\bibitem{alon2004probabilistic}
N.~Alon and J.~H. Spencer, {\em The probabilistic method}.
\newblock John Wiley \& Sons, 2004.

\bibitem{meyer2000matrix}
C.~D. Meyer, {\em Matrix analysis and applied linear algebra}, vol.~71.
\newblock Siam, 2000.

\end{thebibliography}

\end{spacing}
\end{document}